\documentstyle[epsbox,12pt,a4]{amsart} 

\numberwithin{equation}{section}

\newtheorem{thm}{Theorem}
\newtheorem{lem}{Lemma}

\newtheorem{cor}{Corollary}
\newtheorem{rem}{Remark}

\begin{document} 
\title[\null]
{
An analogue of the Chowla--Selberg formula 
for several automorphic $L$-functions
}
\maketitle
\begin{center}
Masatoshi Suzuki \\
\smallskip
\end{center}

\begin{abstract}
In this paper, we will give a certain formula for the Riemann zeta function  
that expresses the Riemann zeta function by an infinte series consisting of $K$-Bessel functions. 
Such an infinite series expression can be regarded as an analogue of the Chowla-Selberg formula.  
Roughly speaking, the Chowla-Selberg formula is  the formula 
that expresses the Epstein zeta-function by an infinite series consisting of $K$-Bessel functions. 
In addition, we also give certain analogues of the Chowla-Selberg formula 
for Dirichlet $L$-functions 
and $L$-functions associated with holomorphic cusp forms. 
Moreover, we introduce a two variable function which is analogous to the real analytic Eisenstein series 
and give a certain limit formula for this one. 
Such a limit formula can be regarded as an analogue of Kronecker's limit formula.  
\end{abstract} 

\section{Introduction and the statement of results} 
Let $\zeta(s)$ be the Riemann zeta function. 
It is defined by the Dirichlet series
\begin{equation}\label{101}
\zeta(s)=\sum_{n=1}^\infty n^{-s}
\end{equation}
for ${\rm Re}(s)>1$ and is extended to a meromorphic function on ${\Bbb C}$. 
In this paper, we will give a certain formula for $\zeta(s)$ 
that expresses $\zeta(s)$ by an infinite series consisting of $K$-Bessel functions. 
It can be regarded as an analogue of the {\it Chowla-Selberg formula} 
for the Epstein zeta-function attached to the binary quadratic form. 

The original Chowla-Selberg formula was given in~\cite{CS}. 
Let $Q(m,n)=am^2+bmn+cn^2$ be a positive definite binary quadratic form 
of discriminant $d=b^2-4ac<0$, $a$, $b$ and $c$ are real numbers and $a>0$.  
The Epstein zeta-function $Z_Q(s)$ attached to the binary quadratic form $Q$ 
is defined by 
\begin{equation}\label{102}
Z_{Q}(s) = \frac{1}{2} \sum_{{(m,n)\in{\Bbb Z}^2}\atop{(m,n)\not=(0,0)}} Q(m,n)^{-s}
\end{equation}
for ${\rm Re}(s)>1$. 
It can be continued meromorphically to the whole complex plane. 
Concerning $Z_Q(s)$, Chowla and Selberg~\cite{CS} gave an identity 
which is called the {\it Chowla-Selberg formula}:
\begin{equation}\label{103}
\aligned
Z_Q(s) 
&= a^{-s} \, \zeta(2s) 
+ a^{-s} \sqrt{\pi} \, \frac{\Gamma(s-\frac{1}{2})}{\Gamma(s)} \, \zeta(2s-1) \, k^{1-2s} + R_Q(s), \\
R_Q(s)
&= \frac{4 a^{-s} k^{-s+\frac{1}{2}}}{\pi^{-s} \Gamma(s)} 
\sum_{n=1}^{\infty} n^{s-\frac{1}{2}} \Bigl( \sum_{d|n} d^{1-2s} \Bigr) 
K_{s-1/2}(2 \pi n k) \cos \bigl(\frac{n \pi b}{a}\bigr),  
\endaligned
\end{equation}
where $K_\nu(z)$ is the $K$-Bessel function (modified Bessel function of the third kind) 
and $k=\sqrt{|d|}/(2a)$. 
The series $R_Q(s)$ converges very quickly. 

The Chowla-Selberg formula has many applications in number theory. 
In particular, 
it can be used to investigate the distribution of the complex zeros of $Z_Q(s)$. 
For example, studies of Stark~\cite{St}, Fujii~\cite{Fu} and Ki~\cite{Ki}. 
Therefore, 
it is expected that if one obtains some suitable analogue 
of the Chowla-Selberg formula for $\zeta(s)$, 
then it can be used to investigate the complex zeros of $\zeta(s)$. 
The first aim of this paper is 
to give an analogue of the Chowla-Selberg formula for $\zeta(s)$ 
motivated by the expectation mentioned above. 

Let $\theta(x)$ be the theta series
\begin{equation}\label{104}
\theta(x)=\sum_{n = -\infty}^{\infty} e^{-\pi n^2 x}. 
\end{equation}
It is well known that the completed Riemann zeta-function 
$\zeta^\ast(s)=\pi^{-s/2}\Gamma(s/2)\zeta(s)$ 
has the integral representation 
\begin{equation}\label{105}
\zeta^\ast(s) = \int_{0}^{\infty} (\theta(t^2) -1) \, t^s \frac{dt}{t}. 
\end{equation}
This holds for ${\rm Re}(s) >1$. 
For a positive real number $\alpha$, 
we define the numbers $c_\alpha(m)$ $(m=1,2,3,\cdots)$ by
\begin{equation}\label{106}
\aligned
\theta(x)^{\alpha}
&= \bigl( 1 + 2 \sum_{n=1}^{\infty} e^{-\pi n^2 x} \bigr)^{\alpha}
= \sum_{j=0}^{\infty} \binom{\alpha}{j} \, 2^j \,  
\bigl(\, \sum_{n=1}^{\infty} e^{-\pi n^2 x} \,\bigr)^j \\
&= 1 + \sum_{m=1}^{\infty} c_{\alpha}(m) e^{-\pi m x}.
\endaligned
\end{equation}
In addition, we define
\begin{equation}\label{107}
Z_\alpha(s)=\int_{0}^{\infty} 
\left( 
\theta(t^2)^\alpha -1 
\right) t^s \frac{dt}{t}.
\end{equation}
This $Z_{\alpha}(s)$ is equal to
the two-variable zeta function $Z_{\Bbb Q}(\alpha, s)$ 
studied by Lagarias and Rains~\cite{LaRa}. 
As a part of their results, the integral on the right-hand side 
converges absolutely for ${\rm Re}(s) > \alpha$ 
and can be continued meromorphically to the whole complex plane. 
By using this $Z_{\alpha}(s)$, our first result is stated as follows. 


\begin{thm}
Let $\alpha$ and $\beta$ be positive real numbers with $\alpha+\beta=1$.  
Then, we have
\begin{equation}\label{thm1:01}
\zeta^\ast(s) 
=  Z_\beta(s)+Z_\alpha(1-s) + \sum_{m=1}^{\infty}\sum_{n=1}^{\infty}
c_{\alpha}(m) c_{\beta}(n) \bigl( \frac{n}{m} \bigr)^{\frac{s-\beta}{4}}
K_{\frac{s-\beta}{2}}\bigl( 2 \pi \sqrt{mn} \,\bigr) \qquad
\end{equation} 
for any $s \in {\Bbb C}$. 
In particular, we have
\begin{equation}\label{thm1:02}
\zeta^\ast(s) 
=  Z_{1/2}(s)+Z_{1/2}(1-s)  +  \sum_{m=1}^{\infty}\sum_{n=1}^{\infty}
c_{1/2}(m) c_{1/2}(n) \bigl( \frac{n}{m} \bigr)^{\frac{s}{4}-\frac{1}{8}}
K_{\frac{s}{2}-\frac{1}{4}}\bigl( 2 \pi \sqrt{mn} \,\bigr) \qquad 
\end{equation}
for any $s \in {\Bbb C}$ by taking $\alpha=\beta=1/2$. 
\end{thm}
Now we compare formula $\eqref{thm1:01}$ in Theorem 1 with 
the Chowla-Selberg formula $\eqref{103}$. 
Let $Q_1(m,n)=m^2+n^2$. Then we have   
\begin{equation}\label{109}
\pi^{-s} \Gamma(s) \, Z_{Q_1}(s) 
= \zeta^\ast(2s) + \zeta^\ast(2-2s) 
 + 4 \sum_{n=1}^{\infty} n^{s-\frac{1}{2}} \Bigl( \sum_{d|n} d^{1-2s} \Bigr) K_{s-1/2}(2 \pi n)
\end{equation}
by using the functional equation $\zeta^\ast(s)=\zeta^\ast(1-s)$. 
On the other hand $\eqref{thm1:01}$ can be written as
\begin{equation}\label{110}
\aligned
\pi^{-s/2} \Gamma(s/2) \, \zeta(s) 
& = Z_\beta(s) +  Z_\alpha(1-s) \\ 
& + \sum_{n=1}^\infty n^{\frac{s-\beta}{4}} 
\Bigl(\, \sum_{d | n} c_\alpha(d) \, c_\beta\bigl(\frac{n}{d}\bigr) \, d^{-\frac{s-\beta}{2}} \Bigr)
\,K_{\frac{s-\beta}{2}}(2 \pi \sqrt{n}\,). 
\endaligned
\end{equation}
Here we observe that formulas $\eqref{109}$ and $\eqref{110}$ have similar forms. 
In this sense, formula $\eqref{thm1:01}$ 
is regarded as an analogue of $\eqref{103}$.  

Also, we can consider formula $\eqref{thm1:01}$ as an 
``one dimensional'' analogue of the Chowla-Selberg formula. 
Formula $\eqref{103}$ is the identity for the Epstein zeta function 
of a two dimensional lattice. 
On the other hand, 
we can regard $\zeta(s)$ as the Epstein zeta function of the one dimensional lattice ${\Bbb Z}$ 
by the formula 
$$
\zeta(s)= \frac{1}{2} \sum_{0 \not= n \in {\Bbb Z}} |n|^{-s}.
$$
In this sense, our $\eqref{thm1:01}$ can be regarded as an one dimensional analogue of $\eqref{103}$. 
The higher dimensional analogue of $\eqref{103}$ had been studied by several authors, 
for example, Terras~\cite{Te}.  
However, as far as the author knows, 
the one dimensional analogue has not been published. 
\smallskip

Before describing the application of Theorem 1, we state two further results. 
Roughly speaking, Theorem 1 is led from only two properties of $\zeta(s)$. 
That two properties are the functional equation $\theta(t^{-2}) = t \, \theta(t^2)$ 
and the integral representation $\eqref{105}$. 
Therefore, if a Dirichlet series $L(s)$ has properties as these, 
then we can prove a formula that is similar to $\eqref{thm1:01}$ for $L(s)$. 
As such examples, we give the following Theorem 2 and Theorem 3 
that are the formulas for Dirichlet $L$-functions 
and $L$-functions attached to holomorphic cusp forms, respectively.    
\smallskip

Let $\chi$ be a real primitive Dirichlet character modulo $q$. 
The theta series $\theta(x,\chi)$ is defined by
\begin{equation}
\theta(x,\chi)=\sum_{n \in {\Bbb Z}} n^{\delta}\chi(n) e^{-\frac{\pi n^2 x}{q}},
\end{equation}
where
\begin{equation}
\delta=
\begin{cases}
0 & \text{if $\chi(-1)=1$,} \\
1 & \text{if $\chi(-1)=-1$.}
\end{cases}
\end{equation}
The completed Dirichlet $L$-function $L^\ast(s,\chi)$ is defined by 
\begin{equation}
L^\ast(s,\chi) = 
q^{\frac{s+\delta}{2}} \pi^{-\frac{s+\delta}{2}} \Gamma\bigl(\frac{s+\delta}{2}\bigr) 
\sum_{n=1}^{\infty}\chi(n) \, n^{-s}.
\end{equation}
The series on the right-hand side converges absolutely for ${\rm Re}(s)>1$.  
It is known that $L^\ast(s,\chi)$ can be continued to an entire function. 
Hereafter, we assume that $\theta(x,\chi)>0$ for any $x>0$. 
Then, for any positive real number $\alpha$, 
we can define the numbers $c_{\chi,\alpha}(m)$ $(m=0,1,2,\cdots)$ by
\begin{equation}
\aligned
\theta(x,\chi)^{\alpha}
&=2^{\alpha}e^{-\frac{\pi \alpha x}{q}}
\bigl(
1 + \sum_{n=2}^{\infty} n^{\delta}\chi(n)e^{-\frac{\pi (n^2-1) x}{q}}  
\bigr)^{\alpha} \\
&=2^{\alpha}e^{-\frac{\pi \alpha x}{q}}
\sum_{j=0}^{\infty} \binom{\alpha}{j}
\bigl(\, \sum_{n=2}^{\infty} n^{\delta}\chi(n)e^{-\frac{\pi (n^2-1) x}{q}} \,\bigr)^j \\
&=2^{\alpha}
\sum_{m=0}^{\infty} c_{\chi,\alpha}(m)e^{-\frac{\pi (m+\alpha) x}{q}}.
\endaligned
\end{equation}


\begin{thm}
Let $\alpha$ and $\beta$ be positive real numbers with $\alpha+\beta=1$.  
Let $\chi$ be a real primitive Dirichlet character modulo $q$. 
Assume that $\theta(x,\chi) > 0$ for any $x>0$. 
Then, we have
\begin{equation}\label{thm2:01}
\aligned
L^\ast & (s,\chi) = 
2  \sum_{m=0}^{\infty}\sum_{n=0}^{\infty}
c_{\chi,\alpha}(m) c_{\chi,\beta}(n) \bigl(\, \frac{n+\beta}{m+\alpha} \,\bigr)^{\frac{s-\beta+\delta(\alpha-\beta)}{4}} \\ 
& \qquad \qquad \qquad \quad \times 
K_{\frac{s-\beta+\delta(\alpha-\beta)}{2}}\bigl(\, \frac{2 \pi}{q} \sqrt{(m+\alpha)(n+\beta)} \,\,\bigr) 
\endaligned
\end{equation}
for any $s \in {\Bbb C}$. 
In particular, we have
\begin{equation}\label{thm2:02}
\aligned
L^\ast & (s,\chi) = 
2  \sum_{m=0}^{\infty}\sum_{n=0}^{\infty}
c_{\chi,1/2}(m) c_{\chi,1/2}(n) \bigl(\, \frac{2n+1}{2m+1} \,\bigr)^{\frac{s}{4}-\frac{1}{8}} 
K_{\frac{s}{2}-\frac{1}{4}}\bigl(\, \frac{\pi}{q} \sqrt{(2m+1)(2n+1)} \,\,\bigr) \qquad 
\endaligned
\end{equation}
by taking $\alpha=\beta=1/2$. 
\end{thm}
\begin{rem}
We need the assumption that $\theta(x,\chi)>0$ for any $x>0$ 
to define $\theta(x,\chi)^{\alpha}=\exp (\alpha \log \theta(x,\chi))$, 
because $\theta(x,\chi)$ is not necessarily positive for all $x>0$ 
for a general real character $\chi$. 
We note that if $\theta(x,\chi)>0$ for any $x>0$ 
then obviously its Mellin transform has no zeros on the real line, 
and hence the corresponding Dirichlet $L$-function has no Siegel zero. 
\end{rem}
By taking suitable binary quadratic form $Q$, 
formula $\eqref{103}$ gives a useful expression of the product $\zeta(s)L(s,\chi)$. 
However it does not give similar type of expressions for each factors $\zeta(s)$ and $L(s,\chi)$. 
In a sense, our formulas $\eqref{thm1:01}$ and $\eqref{thm2:01}$ 
give a factorization of formula $\eqref{103}$. 

Successively, we state the result for $L(s,f)$. 
As for the theory of holomorphic modular forms, 
see Chapter 14 of ~\cite{IK}, for example. 
Let $k$ be a positive integer, 
and let $\Gamma_0(q) = \bigl\{ \bigl( \smallmatrix \ast & \ast \\ c & \ast \endsmallmatrix \bigr) 
\in SL_2({\Bbb Z})
\,|\, c \equiv 0 ~{\rm mod}~q  \bigr\}$. 
Let $f(z)$ be a holomorphic cusp form of weight $k$ on  
$\Gamma_0(q)$, 
and let 
\begin{equation}
f(z)=\sum_{n=1}^{\infty} a_f(n) e^{2\pi i n z}
\end{equation}
be its Fourier expansion at the cusp $i \infty$.
The completed $L$-function $L^\ast(s,f)$ attached to the cusp form $f(z)$ 
is defined by the series
\begin{equation}
L^\ast(s,f)=q^{\frac{s}{2}}(2\pi)^{-s}\Gamma(s)\sum_{n=1}^{\infty}a_f(n) \, n^{-s}. 
\end{equation}
It is known that the series on the right-hand side 
converges absolutely for ${\rm Re}(s)>(k+1)/2$ 
and can be continued to an entire function. 
We also introduce the operator $W$ (sometimes called the Fricke involution) by
\begin{equation}\label{221}
Wf(z)=q^{-\frac{k}{2}}z^{-k}f\bigl(-\frac{1}{qz} \,\bigr). 
\end{equation}
Hereafter, we assume that $a_f(1)=1$ and $f(iy)>0$ for any $y>0$. 
Then, for any positive real number $\alpha$, 
we can define the numbers $c_{f,\alpha}(m)$ $(m=0,1,2,\cdots)$ by
\begin{equation}\label{222}
\aligned
f(iy)^{\alpha}
&=e^{-2\pi \alpha y}
\bigl(
1 + \sum_{n=2}^{\infty} a_f(n) e^{-2\pi (n-1)y}  
\bigr)^{\alpha} \\
& =e^{-2\pi \alpha y}
\sum_{j=0}^{\infty} \binom{\alpha}{j}
\bigl( \, \sum_{n=2}^{\infty} a_f(n)e^{-2\pi (n-1)y} \, \bigr)^j \\
&=
\sum_{m=0}^{\infty} c_{f,\alpha}(m)e^{-2\pi(m+\alpha)y}. 
\endaligned
\end{equation}


\begin{thm}
Let $k$ be a positive integer with $k \equiv 0 \,\, (4)$. 
Let $\alpha$ and $\beta$ be positive real numbers with $\alpha+\beta=1$. 
Let $f(z)$ be a holomorphic cusp form of weight $k$ on $\Gamma_0(q)$, 
and let 
\begin{equation}\label{thm3:01}
f(z)=\sum_{n=1}^{\infty} a_f(n) e^{2\pi i n z}
\end{equation}
be its Fourier expansion at the cusp $i \infty$. 
Assume that $a_f(1)=1$, $Wf=f$ and $f(iy)>0$ for any $y>0$. 
Then, we have 
\begin{equation}\label{thm3:02}
L^\ast(s,f) 
= 2 \, \sum_{m=0}^{\infty}\sum_{n=0}^{\infty}
c_{f,\alpha}(m) c_{f,\beta}(n) \bigl(\, \frac{n+\beta}{m+\alpha} \,\bigr)^{\frac{s-k\beta}{2}}  
K_{s-k\beta}\bigl(\,\frac{4 \pi}{\sqrt{\,q\,}} \sqrt{(m+\alpha)(n+\beta)}\,\,\bigr) \qquad \quad 
\end{equation}
for any $s \in {\Bbb C}$. 
In particular, we have
\begin{equation}\label{thm3:03}
\aligned
L^\ast(s,f) 
&= 2 \, \sum_{m=0}^{\infty}\sum_{n=0}^{\infty}
c_{f,\frac{1}{2}}(m) c_{f,\frac{1}{2}}(n) \bigl(\, \frac{2n+1}{2m+1} \,\bigr)^{\frac{s}{2}-\frac{k}{4}} 
K_{s-\frac{k}{2}}\bigl(\frac{2 \pi}{\sqrt{\,q \,}} \sqrt{(2m+1)(2n+1)}\,\,\bigr) \qquad 
\endaligned
\end{equation}
by taking $\alpha=\beta=1/2$. 
\end{thm}
\begin{rem}
The conditions $k \equiv 0 \,\,(4)$, 
$a_f(1)=1$ and $Wf=f$ are technical ones. 
Even if one of these are not satisfied, 
$\eqref{thm3:02}$ and $\eqref{thm3:03}$ hold up to complex constant multiple. 
See the proof of this theorem in Section $3$. 
On the other hand, the assumptions in Theorem $3$ are actually satisfied 
by several specific cusp forms. 
As an example, the Ramanujan delta function 
$\Delta(z)=e^{2\pi i z} \prod_{n=1}^{\infty} (1-e^{2\pi i n z})^{24} \in S_{12}(\Gamma_0(1))$ 
satisfies the assumptions in this theorem.
\end{rem}

In the latter half of this section, 
we collect three applications and one possible further generalization of the above results. 

\begin{cor} 
Denote $c_{1/2}(n)$ by $c(n)$. 
On the critical line ${\rm Re}(s)=1/2$, 
we have 
\begin{equation*}
\aligned
\zeta^\ast(1/2+it)
= &\, Z_{1/2}(1/2+it) + Z_{1/2}(1/2-it) 
 + \sum_{m=1}^{\infty} c(m)^2 K_{it/2}( 2 \pi m ) \\
& +2 \sum_{1 \leq m<n}^{\infty} c(m)c(n) 
\cos \bigl(\, \frac{t}{4} \, \log \bigl( \frac{n}{m} \bigr) \bigr)
K_{it/2}( 2 \pi \sqrt{mn} \,) \\
= & \, Z_{1/2}(1/2+it) + Z_{1/2}(1/2-it) \\ 
& + \sum_{n=1}^\infty n^{it/4}
\Bigl(\, \sum_{d | n} c(d) \, c(n/d) \, d^{-it/2} \Bigr)
\,K_{it/2}(2 \pi \sqrt{n}\,). 
\endaligned
\end{equation*} 
At the point $t=0$, 
this formula is to be read as the limit $t \to 0$.  
\end{cor}
From the above formula, we have
\begin{equation}\label{fig2}
\aligned
\zeta(1/2+it) & = \zeta_{1/2} (1/2+it)
+ \pi^{1/4}\frac{ \Gamma(it/2)}{\Gamma(1/4+it/2)} \, \zeta_{1/2}(1/2-it) \\
& + \frac{\pi^{1/4+it/2}}{\Gamma(1/4+it/2)} 
\sum_{n=1}^\infty n^{it/4} 
\Bigl(\, \sum_{d|n} c(d) \, c(n/d) \, d^{-it/2} \Bigr)
\,K_{it/2}(2 \pi \sqrt{n}\,), \quad
\endaligned
\end{equation}
where $\zeta_{1/2}(s)=\pi^{s/2} \Gamma(s/2)^{-1} \,Z_{1/2}(s)$. 
As can be seen in Figure 1 and Figure 2 (in the last page), 
the right-hand side gives a good approximation of 
$\zeta(1/2+it)$ on the critical line. 
Our expressions of $\zeta(1/2+it)$ in Corollary 1 or $\eqref{fig2}$ 
may be useful to study $\zeta(s)$ on the critical line.
\begin{cor} 
Let $\chi$ be as in Theorem $2$. Denote $c_{\chi,1/2}(n)$ by $c_\chi(n)$.
The value $L(1,\chi)$ at the edge of the critical strip 
and the central value $L(1/2,\chi)$ have the following infinite series expressions.
\begin{align}
L&(1,\chi)  = \frac{2\pi^{\delta}}{q^{1/2+\delta/2}}
 \sum_{m=0}^{\infty}\sum_{n=0}^{\infty}
c_{\chi}(m) c_{\chi}(n) \bigl(\, \frac{2n+1}{2m+1} \,\bigr)^{\frac{1}{8}} 
K_{\frac{1}{4}}\bigl(\, \frac{\pi}{q} \sqrt{(2m+1)(2n+1)} \,\,\bigr), \qquad \\
L&(1/2, \chi)  = \frac{2\,(\pi/q)^{\frac{1+2\delta}{4}}}{\Gamma(1/4 + \delta/2)\,\,\,\,}
 \sum_{m=0}^{\infty}\sum_{n=0}^{\infty}
c_{\chi}(m) c_{\chi}(n) 
K_{0}\bigl(\, \frac{\pi}{q} \sqrt{(2m+1)(2n+1)} \,\,\bigr).
\end{align}
\end{cor}

\begin{cor} Let $f(z)$ be as in Theorem $3$. Denote $c_{f,1/2}(n)$ by $c_f(n)$.
The value $L((k+1)/2,f)$ at the edge of the critical strip 
has the infinite series expression
\begin{equation}
L((k+1)/2,f) =  \frac{(2\pi)^{(k+1)/2}}{q^{k/4}\Gamma((k+1)/2)} 
 \sum_{m=0}^{\infty}\sum_{n=0}^{\infty}
c_{f}(m) c_{f}(n) \frac{e^{-\frac{2 \pi}{\sqrt{q}} \sqrt{(2m+1)(2n+1)}}}{\sqrt{2m+1}}.
\end{equation}
Here we used $K_{1/2}(x)=\sqrt{\pi/(2x)} \,e^{-x}$. 
Also, the central value $L(k/2,f)$ has the infinite series expression
\begin{equation}
L(k/2,f) = \frac{2}{\Gamma(k/2)} 
\bigl( \frac{2\pi}{\sqrt{\,q\,}} \bigr)^{k/2}
 \sum_{m=0}^{\infty}\sum_{n=0}^{\infty}
c_{f}(m) c_{f}(n) 
K_{0}\bigl(\frac{2 \pi}{\sqrt{\,q \,}} \sqrt{(2m+1)(2n+1)}\,\,\bigr).
\end{equation} 
\end{cor}

The values of $L$-functions at the edge of the critical strip 
and at the central point are very important quantities 
in number theory. 
Our expression for these values may be helpful to investigate these values. 

To generalize Theorem 1 to other directions, 
we recall the real analytic Eisenstein series. 
The real analytic Eisenstein series $E^\ast(z,s)$ is defined by
\begin{equation}
E^\ast(z,s)= \frac{1}{2} \pi^{-s}\Gamma(s) 
\sum_{{(m,n)\in {\Bbb Z}^2}\atop{(m,n)\not=(0,0)}}
\frac{y^s}{|mz+n|^{2s}}
\end{equation}
for $z=x+iy$ with $y>0$ and $s$ with ${\rm Re}(s)>1$. 
It is known that $E^\ast(z,s)$ 
is a modular form of weight zero on ${\rm PSL}_2({\Bbb Z})$ 
as a function of $z$ and has the Fourier expansion
\begin{equation}\label{500-1}
\aligned
E^\ast(z,s)
& = \zeta^\ast(2s) y^s + \zeta^\ast(2-2s) y^{1-s} \\
& + 4 \sqrt{y} \sum_{n=1}^{\infty} n^{s-1/2} \Bigl( \sum_{d|n} d^{1-2s} \Bigr) K_{s-1/2}(2 \pi n y) \cos(2 \pi n x). 
\endaligned
\end{equation}
The Chowla-Selberg formula $\eqref{103}$ is the special case of this Fourier expansion 
at the CM-point $z_Q=\frac{b+\sqrt{d}}{2a}$. 
From this point of view, 
it is natural that 
we regard $\eqref{thm1:01}$ or $\eqref{thm1:02}$ in Theorem 1 as a special case 
of a certain general formula for some two variable function 
analogous to $E^\ast(z,s)$. 
As a candidate of such two variable function, 
we consider the two variable function
\begin{equation}\label{500}
E_{1/2}(y,s)=
\int_{0}^{\infty} \bigl(\sqrt{\theta(y^2t^2)\theta(y^{-2}t^2)} -1 \bigr)\, t^s \frac{dt}{t}.
\end{equation}
This definition resembles the integral representation of $E^\ast(z,s)$ given by
$$
E^\ast(z,s)= \frac{1}{2} \int_{0}^{\infty} \bigl( \Theta_z(t) -1 \bigr) t^s \frac{dt}{t}, 
\quad
\Theta_z(t) = \sum_{(m,n) \in {\Bbb Z}^2} \exp\bigl(-\pi \frac{|mz+n|^2}{y} t \bigr).
$$


\begin{thm}
Let $E_{1/2}(y,s)$ be the function defined by $\eqref{500}$. 
It is defined for $y>0$ and $s$ with ${\rm Re}(s)>1$. 
For any fixed $y>0$, 
$E_{1/2}(y,s)$ can be continued meromorphically to the whole $s$-plane. 
Moreover, the identity
\begin{equation}\label{thm4:01}
\aligned
E_{1/2}(y,s)
& = \, Z_{1/2}(s) \, y^{s} + Z_{1/2}(1-s) \, y^{1-s} \\
& \quad + \sqrt{y} \sum_{n=1}^\infty n^{\frac{s}{4}-\frac{1}{8}} 
\Bigl(\, \sum_{d | n} c_{1/2}(d) \, c_{1/2}(n/d) \, d^{-\frac{s}{2}+\frac{1}{4}} \Bigr)
K_{\frac{s}{2}-\frac{1}{4}} ( 2 \pi n y^2 )
\endaligned
\end{equation}
holds for any $y>0$ and $s \in {\Bbb C}$. 
\end{thm}
At the point $y=1$, we have $E_{1/2}(1,s)=\zeta^\ast(s)$. 
Hence, Theorem 4 contains formula $\eqref{thm1:02}$.  
We can regard the relation between $\eqref{thm1:02}$ and $\eqref{thm4:01}$ 
as an analogue of the relation between $\eqref{109}$ and $\eqref{500-1}$. 
As a corollary of Theorem 4, we obtain the following analogue of Kronecker's limit formula. 
\begin{cor}
For any fixed $y>0$, we have
\begin{equation}
\aligned
\lim_{s \to 1}  \Bigl(\, E_{1/2}(y,s) &- \frac{1}{s-1} \,\Bigr) 
= Z_{1/2}(1) \, y + \gamma_{1/2} \\
& + \sqrt{y} \sum_{n=1}^\infty n^{1/8} 
\Bigl(\, \sum_{d | n} c_{1/2}(d) \, c_{1/2}(n/d) \, d^{-1/4} \Bigr)
K_{1/4} ( 2 \pi n y^2 ),
\endaligned
\end{equation} 
where the constant $\gamma_{1/2}$ is given by 
$\gamma_{1/2}=\underset{s \to 1/2}{\lim} \{ Z_{1/2}(s) - (s - 1/2)^{-1} \}$.
\end{cor}
Original Kronecker's limit formula is the formula 
\begin{equation*}\label{401}
\aligned
\lim_{s \to 1} \Bigl( 2E^\ast(z,s) - \frac{1}{s-1} \Bigr)
= & \,(\, \frac{\pi}{3} - \log 4\pi \,)\, y + \gamma  \\
+ & \, 8 \sqrt{y} \sum_{n=1}^{\infty} \sqrt{n} \, \Bigl( \sum_{d|n} d^{-1} \Bigr) K_{1/2}(2\pi n y)\cos(2 \pi n x),
\endaligned
\end{equation*}
where $\gamma=0.577215\dots$ is Euler's constant. 
It is well known that 
\begin{equation*}\label{402}
8 \sqrt{y} \sum_{n=1}^{\infty} \sqrt{n} \, \Bigl( \sum_{d|n} d^{-1} \Bigr) K_{1/2}(2\pi n y)\cos(2 \pi n x)
= - \frac{\pi}{3} \, y - 4 \log |\eta(z)|,
\end{equation*}
where $\eta(z)$ is the Dedekind eta function
$\eta(z) = e^{\frac{\pi i z}{12}} \prod_{n=1}^{\infty} (1-e^{2\pi i nz})$. 
It seems that it is an interesting problem to consider what is the analogue 
of $\log |\eta(z)|$ (cf. Asai~\cite{Asa}). 


Unfortunately, we have not yet obtained any result 
about application to the distribution of zeros of zeta-functions, 
although it is our first motivation. 
We would like to deal with such applications in our future study. 

The paper is organized as follows. 
In Section $2$, we collect several lemmas and its proofs.  
In Section $3$, we give the proofs of Theorem 1, Theorem 2 and Theorem 3. 
In Section $4$, we prove Theorem 4 in a more general form (Theorem 5).   
\smallskip

\noindent
{\bf Acknowledgement}. 
The author thanks the referee for his detailed and helpful comments on the first version of this paper.

\section{Lemmas}



\begin{lem}
Let $\alpha$ be a positive real number. 
Let $Z_\alpha(s)$ be the function defined by the integral $\eqref{107}$. 
Then, we have
\begin{enumerate} 
\item[(a)] $Z_1(s)=\zeta^\ast(s)$, 
\item[(b)] the integral in $\eqref{107}$ converges absolutely for ${\rm Re}(s) > \alpha$,  
\item[(c)] $Z_\alpha(s)$ is meromorphically continued to the whole $s$-plane, 
\item[(d)] $Z_\alpha(s)$ satisfies the functional equation
\begin{equation}\label{108}
Z_\alpha(s)=Z_\alpha(\alpha -s),
\end{equation}
\item[(e)] in any vertical strip with finite width, 
\begin{equation}
|Z_\alpha(\sigma+it)|=O(t^{-2}) \quad  \text{as $|t| \to \infty$},
\end{equation}
\item[(f)] $Z_\alpha(s)$ is holomorphic except for two simple poles at $s=0$ and $s=\alpha$
with residues $-1$ and $1$, respectively.  
\end{enumerate}
\end{lem}
As mentioned in Section $1$, all properties of $Z_{\alpha}(s)$ in Lemma 1 
have been given in Lagarias and Rains~\cite{LaRa}. 
However, we review the proof of Lemma 1 according to ~\cite{LaRa}, 
because we use a method similar to that to prove succeeding Lemma 2 and Lemma 3. 
\smallskip


\noindent
{\bf Proof of Lemma 1}. 
Assertion (a) follows from definition $\eqref{107}$ 
and the integral representation $\eqref{105}$ of $\zeta(s)$. 
Assertion (b) is an immediate consequence of 
\begin{equation}
\theta(t^2)^{\alpha} -1
=
\begin{cases}\label{306a}
t^{-\alpha} -1 + O(t^N) & \text{as $t \to +0$ for any $N \geq 1$}, \\
O(t^{-N}) & \text{as $t \to +\infty$ for any $N \geq 1$}.
\end{cases}
\end{equation}
Therefore, we prove $\eqref{306a}$. 
We have
\begin{equation}\label{307}
\theta(x)
=\frac{1}{\sqrt{x}} \Bigl(1+2\sum_{n=1}^{\infty} e^{-\pi n^2/x} \Bigr)
\end{equation}
by using the functional equation
\begin{equation}\label{306}
\theta(x^{-1})=\sqrt{x}\,\theta(x). 
\end{equation}
Hence
\begin{equation}\label{308}
\aligned
\theta(t^2)^{\alpha} -1
&= t^{-\alpha} \Bigl( 1 + 2 \sum_{n=1}^\infty e^{-\pi n^2 / t^2} \Bigr)^\alpha  -1 \\
&= t^{-\alpha} \Bigl\{ 
1 + \sum_{j=1}^{\infty} \binom{\alpha}{j} \, 2^j \,  
\bigl(\, \sum_{n=1}^{\infty} e^{-\pi n^2 / t^2} \,\bigr)^j 
\Bigr\} - 1 \\
&= t^{-\alpha} -1 +O(t^N) \quad  \text{as $t \to +0$}
\endaligned
\end{equation} 
for any $N \geq 1$. 
On the other hand, we have 
\begin{equation}\label{308_b}
\aligned
\theta(t^2)^\alpha 
& = \bigl( 1 + 2 \sum_{n=1}^\infty e^{-\pi n^2 t^2} \bigr)^\alpha 
  = 1 + \sum_{j=1}^{\infty} \binom{\alpha}{j} \, 2^j \,  
\bigl(\, \sum_{n=1}^{\infty} e^{-\pi n^2 t^2} \,\bigr)^j \\ 
& = 1 + O(t^{-N}) \quad \text{as $t \to +\infty$}
\endaligned
\end{equation}
for any $N \geq 1$. 
Hence we obtain $\eqref{306a}$ and (b). 
Now we prove (c), (d) and (f) simultaneously. 
Initially,  we suppose that ${\rm Re}(s)>\alpha$ for the convergence of integral. 
We split integral $\eqref{107}$ into 
the part from $1$ to $\infty$ and the part from $0$ to $1$. 
By using the functional equation $\eqref{306}$, we have 
\begin{equation}
\theta(t^{-2})^\alpha -1 = t^{\alpha} (\theta(t^2)^{\alpha} -1) + t^{\alpha}-1.  
\end{equation}
Hence, for the part from $0$ to $1$, we have
\begin{equation*}
\aligned
\int_{0}^{1} (\theta(t^{2})^\alpha -1)\, t^{s} \frac{dt}{t}
&=\int_{1}^{\infty} (\theta(t^{-2})^\alpha -1)\, t^{-s} \frac{dt}{t} \\
&=\int_{1}^{\infty} (\theta(t^{2})^\alpha -1)\, t^{\alpha-s} \frac{dt}{t}
 + \int_{1}^{\infty} (t^{\alpha-s} - t^{-s}) \frac{dt}{t} \\
&=\int_{1}^{\infty} (\theta(t^{2})^\alpha -1)\, t^{\alpha-s} \frac{dt}{t}
+\frac{1}{s-\alpha}-\frac{1}{s}.
\endaligned
\end{equation*}
This leads 
\begin{equation}\label{311}
Z_\alpha(s)=\frac{1}{s-\alpha}-\frac{1}{s}
 + \int_{1}^{\infty} (\theta(t^{2})^\alpha -1)\, (t^s+t^{\alpha-s}) \frac{dt}{t}.
\end{equation}
Here, we find that the integral on the right-hand side converges absolutely for any $s \in {\Bbb C}$ 
by $\eqref{306a}$. 
Therefore the integral is an entire function on ${\Bbb C}$. 
Hence, the representation $\eqref{311}$ of $Z_\alpha(s)$ 
gives the meromorphic continuation of $Z_\alpha(s)$ to the whole complex plane. 
More precisely, $Z_\alpha(s)$ is holomorphic 
except for simple poles at $s=0$ and $s=\alpha$ 
with residues $-1$ and $1$ respectively.  
Because the right-hand side of $\eqref{311}$ is invariant under $s \mapsto \alpha - s$, 
we obtain the functional equation $\eqref{108}$. We complete the proof of (c), (d) and (f). 
Finally, we prove (e). 
For $s$ with ${\rm Re}(s) \not=0,\alpha$, we have
\begin{equation}\label{est_ver_1}
\frac{1}{s-\alpha} - \frac{1}{s} = O\bigl( \frac{1}{1+|{\rm Im}(s)|^2} \bigr).
\end{equation}
Also, we have
\begin{equation}\label{est_ver_2}
 \int_{1}^{\infty} (\theta(t^{2})^\alpha -1)\, (t^s+t^{\alpha-s}) \frac{dt}{t}
 = O((1+|{\rm Im}(s)|)^{-N})
\end{equation}
for any $N \geq 1$  
in any vertical strip with finite width 
by using integration by parts. 
Here, we used 
\begin{equation}\label{327}
\frac{d^k}{dt^k}( \theta(t^2)^\alpha - 1)=
O(t^{-N}) \quad  \text{as $t \to +\infty$} 
\end{equation} 
for any $N \geq 1$ and $k \geq 0$ which is obtained by $\eqref{308_b}$. 
These two estimates and  $\eqref{311}$ imply assertion (e). 
\hfill $\Box$
\bigskip


\begin{lem} 
Let $\alpha$ be a positive real number. 
Let $\chi$ be a real primitive Dirichlet character modulo $q$. 
Assume that $\theta(x,\chi) > 0$ for any $x>0$. 
Define
\begin{equation}\label{lem2:01}
Z_\alpha(s,\chi)=\int_{0}^{\infty}  
\theta(t^2,\chi)^\alpha 
\, t^{s+ \delta \alpha} \frac{dt}{t}.
\end{equation}
Then, we have
\begin{enumerate}
\item[(a)] $Z_1(s,\chi)=L^\ast(s,\chi)$, 
\item[(b)] the integral in $\eqref{lem2:01}$ converges absolutely for any $s \in {\Bbb C}$,
\item[(c)] $Z_\alpha(s,\chi)$ is an entire function, 
\item[(d)] $Z_\alpha(s,\chi)$ satisfies the functional equation 
\begin{equation}\label{lem2:02}
Z_\alpha(s,\chi) =  Z_\alpha(\alpha -s,\chi),
\end{equation} 
\item[(e)] in any vertical strip with finite width,
\begin{equation}
|Z_\alpha(\sigma+it,\chi)|=O(|t|^{-N}) 
\quad \text{as $|t| \to +\infty$}
\end{equation}
for any $N \geq 1$. 
\end{enumerate}
\end{lem}


\noindent
{\bf Proof of Lemma 2}. 
First, we note that 
$\theta(x,\chi)^\alpha = \exp( \alpha \log \theta(x,\chi))$ is well defined by 
the assumption that $\theta(x,\chi)>0$ for any $x>0$. 
Assertion (a) is a consequence of the integral representation
\begin{equation}\label{322a}
L^\ast(s,\chi)=\int_{0}^{\infty} \theta(t^2,\chi) \, t^{s+\delta} \, \frac{dt}{t}
\end{equation}
which can be found in~\cite[\S 4.6]{IK}, for example.  
To prove (b), we show that 
\begin{equation}\label{324}
\theta(x,\chi)=
\begin{cases}
O(x^{N}) & \text{as $x \to +0$ for any $N \geq 1$}, \\
O(x^{-N}) & \text{as $x \to +\infty$ for any $N \geq 1$.}
\end{cases}
\end{equation} 
Clearly, this implies (b), since $\alpha$ in $\eqref{lem2:01}$ is positive. 
From the definition, we have 
\begin{equation}\label{322_a}
\theta(x,\chi)=2\sum_{n=1}^\infty n^\delta \chi(n) e^{-\frac{\pi n^2}{q} x } 
= O(x^{-N}) \quad \text{as $x \to +\infty$}
\end{equation}
for any $N \geq 1$. 
On the other hand, we have 
\begin{equation}\label{323}
\theta(x,\chi)
= 2 \, x^{-\frac{1}{2}-\delta} \sum_{n=1}^{\infty} n^{\delta} \chi(n)  e^{-\frac{\pi n^2}{q x}}
= O(x^N) \quad  \text{as $x \to +0$}
\end{equation}
for any $N \geq 1$ 
by using the functional equation 
\begin{equation}\label{322}
\theta(x^{-1},\chi) = x^{\frac{1}{2}+\delta}\,\theta(x,\chi). 
\end{equation}
In $\eqref{322}$, we used the fact that the Gauss sum $\tau(\chi)$ is equal to $i^{\delta}\sqrt{\,q \,}$ 
for any real primitive character $\chi$ modulo $q$. 
For this fact and the functional equation of $\theta(x,\chi)$, see ~\cite[\S 4.6]{IK}. 
From $\eqref{322_a}$ and $\eqref{323}$, we obtain $\eqref{324}$.  
Assertion (c) is an immediate consequence of (b). 
To prove (d), we split integral $\eqref{lem2:01}$ into 
the part from $1$ to $\infty$ and the part from $0$ to $1$. 
By using the functional equation $\eqref{322}$, we have
\begin{equation}\label{325}
\aligned
\int_{0}^{1} \theta(t^{2},\chi)^\alpha \, t^{s + \delta\alpha} \frac{dt}{t}
= \int_{1}^{\infty} \theta(t^{-2},\chi)^\alpha \, t^{-s-\delta\alpha} \frac{dt}{t} 
= \int_{1}^{\infty} \theta(t^{2},\chi)^\alpha \, t^{\alpha-s+\delta\alpha} \frac{dt}{t}. \quad
\endaligned
\end{equation}
Therefore, 
\begin{equation}\label{326}
Z_\alpha(s,\chi) = 
\int_{1}^{\infty} \theta(t^{2},\chi)^\alpha \, t^{s+\delta\alpha} \frac{dt}{t}
+ \int_{1}^{\infty} \theta(t^{2},\chi)^\alpha \, t^{\alpha-s +\delta\alpha} \frac{dt}{t}.
\end{equation} 
The integrals on the right-hand side converge absolutely for any $s \in {\Bbb C}$, 
because of $\eqref{324}$. 
Since the right-hand side of $\eqref{326}$ is invariant when we replace $s$ by $\alpha - s$, 
we obtain the functional equation $\eqref{lem2:02}$. 
Finally, we prove (d). 
From the definition of $\theta(x,\chi)$, we have 
\begin{equation}
\aligned
\theta(t^2,\chi)^\alpha 
& = 2^\alpha e^{- \frac{\pi \alpha}{q} t^2 } 
\Bigl( 1 +  \sum_{n=2}^\infty n^{\delta} \chi(n) e^{-\frac{\pi (n^2-1)}{q} t^2} \Bigr)^\alpha \\
& = 2^\alpha e^{- \frac{\pi \alpha}{q} t^2 } 
\Bigl\{ 1 + \sum_{j=1}^{\infty} \binom{\alpha}{j} \,   
\bigl(\, \sum_{n=2}^{\infty} n^{\delta} \chi(n) e^{-\frac{\pi (n^2-1)}{q} t^2} \,\bigr)^j 
\Bigr\}.
\endaligned
\end{equation}
This shows that 
\begin{equation}\label{327}
\frac{d^k}{dt^k} \theta(t^2,\chi)^\alpha =
O(t^{-N}) \quad \text{as $t \to +\infty$} 
\end{equation} 
for any $N \geq 1$ and $k \geq 0$. 
Hence, we obtain (d) from $\eqref{326}$ by using $\eqref{327}$ and integration by parts. 
\hfill $\Box$
\bigskip


\begin{lem}  
Let $k$ be a positive integer with $k \equiv 0 \,\, (4)$ 
and let $f(z)$ be a holomorphic cusp form of weight $k$ on $\Gamma_0(q)$. 
Assume that $f(iy)>0$ for any $y>0$. 
Let $\alpha$ be a positive real number. 
Define
\begin{equation}\label{lem3:02}
Z_\alpha(s,f)=\int_{0}^{\infty}  
f\bigl(\, \frac{iy}{\sqrt{q}} \,\bigr)^\alpha \, y^{s} \frac{dy}{y}.
\end{equation}
Then we have
\begin{enumerate}
\item[(a)] $Z_1(s,f)=L^\ast(s,f)$, 
\item[(b)] the integral in $\eqref{lem3:02}$  converges absolutely for any $s \in {\Bbb C}$, 
\item[(c)] $Z_\alpha(s,f)$ is an entire function, 
\item[(d)] $Z_\alpha(s,f)$ satisfies the functional equation 
\begin{equation}\label{lem3:03}
Z_\alpha(s,f) = Z_\alpha(k\alpha -s,Wf).
\end{equation}
\item[(e)] in any vertical strip with finite width,
\begin{equation}
|Z_\alpha(\sigma+it,f)|=O(|t|^{-N}) 
\quad \text{as $|t| \to +\infty$}
\end{equation}
for any $N \geq 1$. 
\end{enumerate}
\end{lem}


\noindent
{\bf Proof of Lemma 3}. 
First, we note that 
$f(\frac{iy}{\sqrt{q}})^\alpha = \exp( \alpha \log f(\frac{iy}{\sqrt{q}}))$ is well defined by 
the assumption that $f(iy)>0$ for any $y>0$. 
Assertion (a) is a consequence of the integral representation
\begin{equation}
L^\ast(s,f)=\int_{0}^{\infty} f\bigl(\, \frac{iy}{\sqrt{q}} \,\bigr) y^s \frac{dy}{y}
\end{equation}
which  can be found in~\cite[\S 14.5]{IK}, for example. 
It is known that if $f$ belongs to $S_k(\Gamma_0(q))$, then $Wf$ also belongs to $S_k(\Gamma_0(q))$. 
Therefore both $f$ and $Wf$ are cusp forms, i.e., 
$f(iy)$ and $Wf(iy)$ decay exponentially fast as $y \to +\infty$.    
Hence equation $\eqref{221}$ shows that 
integral $\eqref{lem3:02}$ converges absolutely for any $s \in {\Bbb C}$. 
This is assertion (b). 
Assertion (c) is an immediate consequence of (b). 
To prove (d), we split integral $\eqref{lem3:02}$ 
into the part from $1$ to $\infty$ and the part from $0$ to $1$.
For the part from $0$ to $1$, we have 
\begin{equation}\label{334}
\int_{0}^{1} f\bigl(\frac{iy}{\sqrt{\,q \,}}\bigr)^\alpha \, y^{s} \frac{dy}{y}
= \int_{1}^{\infty} f\bigl(\frac{-1}{\sqrt{q}\,iy}\bigr)^\alpha \, y^{-s} \frac{dy}{y} 
= \int_{1}^{\infty} \Bigl(Wf\bigl(\frac{iy}{\sqrt{\,q \,}}\bigr)\Bigr)^\alpha \, y^{k \alpha-s} \frac{dy}{y}
\quad
\end{equation}
by using $\eqref{221}$ and $k \equiv 0 \,\, (4)$ in the second equality.  
Therefore, we obtain
\begin{equation}\label{335}
Z_\alpha(s,f) = 
\int_{1}^{\infty} f\bigl(\frac{iy}{\sqrt{\,q \,}}\bigr)^\alpha \, y^{s} \frac{dy}{y}
+ \int_{1}^{\infty} \Bigl(Wf\bigl(\frac{iy}{\sqrt{\,q \,}}\bigr)\Bigr)^\alpha \, y^{k \alpha-s} \frac{dy}{y}.
\end{equation}
Because $W(Wf)=f$, representation $\eqref{335}$ shows the functional equation $\eqref{lem3:03}$. 
This is assertion (d). 
By a way similar to the proof of Lemma 2, 
we find that 
$\frac{d^k}{dy^k} f(\frac{iy}{\sqrt{q}})^\alpha$ and 
$\frac{d^k}{dy^k} ( Wf( \frac{iy}{\sqrt{q}}) )^\alpha$ 
decay exponentially fast as $y \to +\infty$ for any $k \geq 0$. 
Hence, we obtain (e) from $\eqref{335}$ by using integration by parts. 
\hfill $\Box$
\bigskip


\begin{lem}
Let $f(x)$ and $g(x)$ be continuous functions on $(0,\infty)$. 
Let $F(w)$ and $G(w)$ be the Mellin transforms 
of $f(x)$ and $g(x)$, respectively. 
Suppose that 
$f(x)$ and $g(x)$ 
decay rapidly as $x \to \infty$, 
and $f(x)=O(x^{-\alpha})$ as $x \to +0$, 
$g(x)=O(x^{-\beta})$ as $x \to +0$, 
where $\alpha$ and $\beta$ are real numbers. 
In addition, we suppose that 
$F(c+it)$ belongs to $L^1(-\infty,\infty)$ for any $c>\alpha$ as a function of $t$ . 
Then  
\begin{equation}\label{lem4}
\frac{1}{2 \pi i} \int_{c - i \infty}^{c + i \infty} 
F(w) \, G(s-w) \, dw = \int_{0}^{\infty} f(x) \, g(x) \, x^{s} \frac{dx}{x}
\end{equation}
for ${\rm Re}(s) > \alpha + \beta$, 
where $c$ is chosen as $\alpha < c < {\rm Re}(s) - \beta$. 
\end{lem}

\noindent
{\bf Proof of Lemma 4.} 
The Mellin transform $F(w)$ of $f(x)$ is defined by
$$
F(w) = \int_{0}^{\infty} f(x) x^w \frac{dx}{x}. 
$$
Under the assumption the integral on the right-hand side 
converges absolutely for ${\rm Re}(w) > \alpha$. 
Hence it follows from the Mellin inversion formula that 
$$
\frac{1}{2 \pi i} \int_{c-i \infty}^{c + i \infty} F(w) x^{-w} dw = f(x), 
$$
where $c > \alpha$. 
Multiply $g(x) x^{s-1}$ to the both sides 
and integrate over $(0,\infty)$. 
Then we have formally
\begin{equation}\label{302}
\int_{0}^{\infty} \left(\frac{1}{2 \pi i} \int_{c-i \infty}^{c + i \infty} F(w) x^{-w} dw \right) g(x)x^{s} \frac{dx}{x} 
= \int_{0}^{\infty} f(x)g(x) x^{s} \frac{dx}{x}. 
\end{equation}
This equality is valid for $s$ with ${\rm Re}(s) > \alpha + \beta$. 
In fact the integral on the right-hand side converges absolutely for ${\rm Re}(s)>\alpha+\beta$, 
because $f(x)g(x)=O(x^{-\alpha-\beta})$ as $x \to +0$ and decays rapidly as $x \to +\infty$. 
On the left-hand side of $\eqref{302}$, we have formally
\begin{equation}\label{302a}
\aligned
\int_{0}^{\infty} \left(\frac{1}{2 \pi i} 
\int_{c-i \infty}^{c + i \infty} F(w) x^{-w} dw \right) 
g(x)x^{s} \frac{dx}{x}
&=\frac{1}{2 \pi i} \int_{c - i \infty}^{c + i \infty} F(w)
\left( \int_{0}^{\infty} g(x)x^{s-w} \frac{dx}{x} \right) dw \\
&= \frac{1}{2 \pi i} \int_{c - i \infty}^{c + i \infty} F(w)G(s-w)dw. 
\endaligned
\end{equation}
This is justified by Fubini's theorem if the integral 
\begin{equation}\label{lem4_wint}
\int_{0}^{\infty} \int_{-\infty}^{\infty} F(c+it) \, g(x) x^{s-c-it-1} dt dx  
\end{equation}
converges absolutely. 
From the assumptions for $F(w)$ and $g(x)$, 
the double integral $\eqref{lem4_wint}$ converges absolutely 
when $c$ is chosen as $\alpha < c < {\rm Re}(s)-\beta$. 
Hence we obtain $\eqref{lem4}$ 
for $s$ with ${\rm Re}(s)>\alpha+\beta$ and 
$c$ with $\alpha < c < {\rm Re}(s)-\beta$, 
by  $\eqref{302}$ and $\eqref{302a}$. \hfill $\Box$
\bigskip


\begin{lem}
We have
\begin{equation}\label{lem5:02}
\frac{1}{2\pi i} \int_{c - i \infty}^{c + i \infty} 
\Gamma\bigl(\frac{s}{2}\bigr)\Gamma\bigl(\frac{s-\nu}{2}\bigr) \, x^{-s}ds
= 4 \, K_{\frac{\nu}{2}}( 2 x ) \, x^{-\frac{\nu}{2}}
\end{equation} 
and 
\begin{equation}\label{lem5:03}
\frac{1}{2\pi i} \int_{c - i \infty}^{c + i \infty} 
\Gamma(s)\Gamma(s-\nu) \, x^{-s}ds
= 2 \, K_\nu( 2 \sqrt{x} \,) \, x^{-\frac{\nu}{2}}
\end{equation}
for $c > {\rm max}\{0, {\rm Re}(\nu) \}$, where $K_\nu(x)$ is the $K$-Bessel function 
defined by 
$$
K_\nu(x)=\frac{1}{2} \int_{0}^{\infty} e^{-\frac{x}{2} ( t + \frac{1}{t} )} t^{-\nu} \frac{dt}{t} \quad (x>0). 
$$ 
\end{lem}

\noindent
{\bf Proof of Lemma 5.} 
Equations $\eqref{lem5:02}$ and $\eqref{lem5:03}$ 
are simple modifications of the equation in \cite[p.197, l.10]{Ti}. 
In \cite[p.197, l.10]{Ti}, 
we replace $\nu$ by $\nu/2$ and $x$ by $2x$. 
Then, we obtain
$$
4x^{-\frac{\nu}{2}}K_{\frac{\nu}{2}}(2x)
= \frac{1}{2\pi i} \int_{c - i \infty}^{c + i \infty} 
\Gamma\bigl( \frac{s}{2} \bigr)
\Gamma\bigl( \frac{s-\nu}{2} \bigr) x^{-s} ds
$$
for any real $\nu$, where $c$ is chosen as $c > {\rm max}\{0, \nu \}$. 
On the other hand, in \cite[p.197, l.10]{Ti}, 
we replace $s$ by $2s$ and $x$ by $2\sqrt{x}$. 
Then, we obtain
$$
2 \, K_\nu( 2 \sqrt{x} \,) \, x^{-\frac{\nu}{2}} 
= \frac{1}{2\pi i} \int_{c - i \infty}^{c + i \infty} 
\Gamma(s)\Gamma(s-\nu) \, x^{-s}ds 
$$
for any real $\nu$, where $c$ is chosen as $c > {\rm max}\{0, \nu \}$. 
For any complex $\nu$, we find that 
the integrals on the right-hand side of the above converge absolutely 
for $c>{\rm max}\{0,{\rm Re}(\nu)\}$ 
by using Stirling's formula. 
Hence $\eqref{lem5:02}$ and $\eqref{lem5:03}$ hold for $c>{\rm max}\{0,{\rm Re}(\nu)\}$ 
by analytic continuation. \hfill $\Box$


\section{Proofs of Theorem 1 to Theorem 3}

First, we describe the proof of Theorem 1 in detail. 
Because the latter theorems are proved by a very similar way, 
we only describe an outline for the proof of Theorem 2 and Theorem 3.


\subsection{Proof of Theorem 1} 
The first step of the proof is to prove that the identity
\begin{equation}\label{320}
\zeta^\ast(s)
=Z_\beta(s)+Z_\alpha(s-\beta)+
\frac{1}{2\pi i}\int_{c -i \infty}^{c + i \infty}
Z_{\alpha}(w)Z_\beta(s-w) dw
\end{equation}
holds for positive $\alpha$ and $\beta$ with $\alpha+\beta=1$ under the conditions ${\rm Re}(s)>1$, $c>{\rm Re}(s)$.  
To prove $\eqref{320}$, we calculate the integral
\begin{equation}\label{315}
I(s;\alpha,\beta)=\int_{0}^{\infty} (\theta(t^2)^\alpha -1)(\theta(t^2)^\beta -1)\,t^{s}\frac{dt}{t}
\end{equation}
in two ways. We start with formal computations for $I(s;\alpha,\beta)$. 
Since 
\begin{equation}\label{316}
(\theta(t^2)^\alpha -1)(\theta(t^2)^\beta -1)= 
(\theta(t^2)-1)-(\theta(t^2)^\alpha-1)-(\theta(t^2)^\beta -1),
\end{equation}
we have 
\begin{equation}\label{317}
I(s;\alpha,\beta)=\zeta^\ast(s)-Z_\alpha(s)-Z_\beta(s).
\end{equation}
On the other hand, we have 
\begin{equation}\label{318}
I(s;\alpha,\beta)
=\frac{1}{2\pi i}\int_{c -i \infty}^{c + i \infty}
Z_{\alpha}(w)Z_\beta(s-w) dw 
\end{equation}
by applying Lemma 4 formally to $\theta(t^2)^\alpha-1$ and $\theta(t^2)^\beta-1$,  
where we understand that $c$ is chosen as $\alpha<c<{\rm Re}(s)-\beta$. 
By moving the path of integration to the vertical line 
${\rm Re}(w)=c^\prime > {\rm Re}(s)$, we have 
\begin{equation}\label{319}
\aligned
\frac{1}{2\pi i}\int_{c -i \infty}^{c + i \infty}
 & Z_{\alpha}(w)Z_\beta(s-w) dw \\
= & \, Z_\alpha(s-\beta)-Z_\alpha(s)  +
\frac{1}{2\pi i}\int_{c^\prime -i \infty}^{c^\prime + i \infty}
Z_{\alpha}(w)Z_\beta(s-w) dw.
\endaligned
\end{equation}
By putting $\eqref{317}$, $\eqref{318}$, $\eqref{319}$ together, 
we obtain $\eqref{320}$ in a formal sense. 

Now we justify each steps. 
Because $\alpha$ and $\beta$ are positive and $\alpha+\beta=1$, 
we have
$$
(\theta(t^2)^\alpha -1)(\theta(t^2)^\beta -1) 
=
\begin{cases} 
O(t^{-1}) & \text{as $t \to +0$}, \\
O(t^{-N}) & \text{as $t \to +\infty$ for any $N \geq 1$}. 
\end{cases}
$$ 
Therefore, the integral $I(s;\alpha,\beta)$ is defined for ${\rm Re}(s)>1$. 
From (b) of Lemma 1, $\eqref{317}$ holds for ${\rm Re}(s)>1$. 
To apply Lemma 4 to $f(t)=\theta(t^2)^\alpha -1$ and $g(t)=\theta(t^2)^\beta -1$, 
we check the conditions in Lemma 4 for these $f(t)$ and $g(t)$.  
Because 
$$
\theta(t^2)^\alpha -1
=
\begin{cases} 
O(t^{-\alpha}) & \text{as $t \to +0$}, \\
O(t^{-N}) & \text{as $t \to +\infty$ for any $N \geq 1$},  
\end{cases}
$$ 
our $f(t)$ and $g(t)$ satisfy the growth condition in Lemma 4. 
Also, the Mellin transform $F(w)$ of $f(t)$ is defined for ${\rm Re}(w)>\alpha$, 
and $F(w)$ is holomorphic in ${\rm Re}(w)>\alpha$. 
Therefore, the estimate in (e) of Lemma 1 implies $F(c+it) \in L^1(-\infty,\infty)$ for any $c>\alpha$. 
Hence we can apply Lemma 4 to $f(t)=\theta(t^2)^\alpha -1$ and $g(t)=\theta(t^2)^\beta -1$. 
As the result, we obtain $\eqref{318}$ with ${\rm Re}(s)>1$ and $\alpha < c <{\rm Re}(s) -\beta$.  
From the estimate in (e) of Lemma 1, 
we can move the path of integration in $\eqref{318}$ 
to the vertical line ${\rm Re}(w)=c^\prime > {\rm Re}(s)$. 
From (f) of Lemma 1, as a function of $w$, 
the poles of $Z_\alpha(w)Z_\beta(s-w)$ 
in the vertical strip $c < {\rm Re}(w) < {\rm Re}(s)+1$ 
are only the simple poles $w=s-\beta$ and $w=s$. 
Their residues are $-Z_\alpha(s-\beta)$ and $Z_\alpha(s)$, respectively. 
Hence we obtain $\eqref{319}$ by using the residue theorem. 
Together with the above, we obtain $\eqref{320}$ 
for positive $\alpha$ and $\beta$ with $\alpha+\beta=1$ 
under the conditions ${\rm Re}(s)>1$, $c>{\rm Re}(s)$.
We complete the first step of the proof. 
\smallskip

The Second step of the proof is to express the integral in $\eqref{320}$ 
by the series involving $K$-Bessel functions. 
To do this, we need an estimate for $c_\alpha(n)$ defined by $\eqref{106}$. 
Let $\alpha>0$ and let $\vartheta(\tau)=\sum_{n \in {\Bbb Z}}e^{\pi i n^2 \tau}$ for ${\rm Im}(\tau)>0$. 
Then $\theta(x)=\vartheta(ix)$ and $\vartheta(\tau)^\alpha$ is a modular form of real weight $\alpha/2$ 
on the theta group with a unitary multiplier system. 
A classical estimate of Petersson~\cite{Pe} and Lehner~\cite{Le} for the Fourier coefficients of 
arbitrary modular forms of positive real weight with multiplier system 
show that they grow polynomially in $m$ as
\begin{equation}\label{313}
c_\alpha(m) =
O(m^{\alpha/4})
\end{equation} 
for $0< \alpha < 4$, where the $O$-constant depends on $\alpha$ in an specified manner. 
By using this estimate, we calculate the integral in $\eqref{320}$. 
From $\eqref{313}$, we obtain
\begin{equation}\label{314}
Z_\alpha(s)
=\sum_{m=1}^{\infty} c_{\alpha}(m) \int_{0}^{\infty} e^{-\pi m t^2} t^{s} \frac{dt}{t} 
=\frac{1}{2} \, \pi^{-s/2}\Gamma(s/2)
\sum_{m=1}^{\infty} c_{\alpha}(m) m^{-s/2}
\end{equation}
for ${\rm Re}(s)>(\alpha/2)+2$ with $0 < \alpha < 4$. 
On the other hand, 
the integrand $Z_\alpha(w) Z_{\beta}(s-w)$ 
on the right-hand side of $\eqref{320}$ 
is equal to $Z_\alpha(w)Z_\beta(w-s+\beta)$ 
because of the functional equation $\eqref{108}$. 
Hence we obtain 
\begin{equation}\label{321}
\aligned
Z_\alpha(w) & Z_{\beta}(s-w)
= Z_\alpha(w)Z_\beta(w-s+\beta) \\
= &
\frac{1}{4} \, \sum_{m=1}^{\infty} \sum_{n=1}^{\infty} 
c_{\alpha}(m) c_{\beta}(n) (\pi n)^{\frac{s-\beta}{2}} 
\times 
\Gamma\bigl(\frac{w}{2}\bigr) 
\Gamma\bigl(\frac{w-s+\beta}{2}\bigr)(\pi \sqrt{mn}\,)^{-w}
\endaligned
\end{equation}
when ${\rm Re}(w)$ is sufficiently large. 
If necessary, by moving the path of integration in $\eqref{320}$ to the right, 
we have
\begin{equation}\label{321_1}
\aligned
\frac{1}{2\pi i} & \int_{c -i \infty}^{c + i \infty}
Z_{\alpha}(w)Z_\beta(s-w) dw \\
& =
\frac{1}{4} \, \sum_{m=1}^{\infty} \sum_{n=1}^{\infty} 
c_{\alpha}(m) c_{\beta}(n) (\pi n)^{\frac{s-\beta}{2}}  
\frac{1}{2\pi i} \int_{c -i \infty}^{c + i \infty}
\Gamma\bigl(\frac{w}{2}\bigr) \Gamma\bigl(\frac{w-s+\beta}{2}\bigr)(\pi \sqrt{mn}\,)^{-w} dw \\
& = \sum_{m=1}^{\infty} \sum_{n=1}^{\infty} 
c_{\alpha}(m) c_{\beta}(n) \bigl( \frac{n}{m} \bigr)^{\frac{s-\beta}{4}} 
K_{\frac{s-\beta}{2}}(2\pi \sqrt{mn} \,).
\endaligned
\end{equation}
Here we used identity $\eqref{lem5:02}$ in Lemma 5. 
Now we obtain 
\begin{equation}\label{almost_thm1}
\zeta^\ast(s)
=Z_\beta(s)+Z_\alpha(s-\beta)+
\sum_{m=1}^{\infty} \sum_{n=1}^{\infty} 
c_{\alpha}(m) c_{\beta}(n) \bigl( \frac{n}{m} \bigr)^{\frac{s-\beta}{4}} 
K_{\frac{s-\beta}{2}}(2\pi \sqrt{mn} \,)
\end{equation}
for positive $\alpha$ and $\beta$ with $\alpha+\beta=1$ 
under the condition ${\rm Re}(s)>1$. 
The series on the right-hand side converges absolutely for any $s \in {\Bbb C}$, 
because $c_\alpha(m)$ is of polynomial order and 
the asymptotic behavior of the $K$-Bessel function is given by    
$$
K_\nu(x) = \sqrt{\frac{\pi}{2x}} \, e^{-x}(1+O(x^{-1}))
$$
for any fixed $\nu \in {\Bbb C}$ (see~\cite{Wa}).  
Hence $\eqref{almost_thm1}$ holds 
on the whole complex plane 
for positive $\alpha$ and $\beta$ with $\alpha+\beta=1$ 
by analytic continuation. 
By using the functional equation $\eqref{108}$, 
$Z_\alpha(s-\beta)=Z_{\alpha}(\alpha+\beta-s)=Z_\alpha(1-s)$. 
Therefore, 
$\eqref{almost_thm1}$ is equivalent to identity $\eqref{thm1:01}$ in Theorem 1. 
We complete the proof of Theorem 1. 
\hfill $\Box$


\subsection{Proof of Theorem 2} 
As in the proof of Theorem 1, we first prove that 
\begin{equation}\label{332}
L^\ast(s,\chi)
=\frac{1}{2\pi i}\int_{c -i \infty}^{c + i \infty}
Z_{\alpha}(w, \chi)Z_\beta(s-w, \chi) dw
\end{equation}
holds for any positive $\alpha$ and $\beta$ with $\alpha+\beta=1$ and $s \in {\Bbb C}$, 
where $c$ is an arbitrary real number. 
Unlike the situation of $\zeta(s)$, this holds for any $s \in {\Bbb C}$. 
To prove $\eqref{332}$, 
we define 
\begin{equation}\label{329}
I_\chi(s;\alpha,\beta)=\int_{0}^{\infty} \theta(t^2,\chi)^\alpha \theta(t^2,\chi)^\beta \,
t^{s+\delta(\alpha+\beta)}\frac{dt}{t}
\end{equation}
for positive real numbers $\alpha$ and $\beta$ with $\alpha+\beta=1$.
The integral on the right-hand side converges absolutely for any $s \in {\Bbb C}$ by $\eqref{324}$. 
We have
\begin{equation}\label{330}
I_\chi(s;\alpha,\beta)=L^\ast(s,\chi) 
\end{equation}
for any $s\in {\Bbb C}$, because of $\alpha+\beta=1$ and $\eqref{322a}$. 
While, from (d) of Lemma 2, 
we can apply Lemma 4 to $f(t)=\theta(t^2,\chi)^\alpha$ and $g(t)=\theta(t^2,\chi)^\beta$. 
From $\eqref{324}$, we obtain 
\begin{equation}\label{331}
I_\chi(s;\alpha,\beta)
=\frac{1}{2\pi i}\int_{c -i \infty}^{c + i \infty}
Z_{\alpha}(w,\chi)Z_\beta(s-w,\chi) dw, 
\end{equation}
for 
${\rm Re}(s) > - 2 \alpha N - 2 \beta N$, 
$-2 \alpha N < c < {\rm Re}(s) + 2 \beta N$, 
where $N$ is an arbitrary positive integer. 
Hence $\eqref{331}$ holds for any $s \in {\Bbb C}$ and $c \in {\Bbb R}$. 
From $\eqref{330}$ and $\eqref{331}$, 
we obtain $\eqref{332}$ under the desired conditions.
\smallskip

Now, we calculate the integral on the right-hand side of $\eqref{332}$. 
Let $\vartheta(\tau,\chi)=\sum_{n \in {\Bbb Z}} n^{\delta} \chi(n) e^{\pi i n^2 \tau}$ for ${\rm Im}(\tau)>0$. 
Then $\theta(x,\chi)=\vartheta(ix,\chi)$ and $\vartheta(\tau,\chi)^\alpha$ is a modular form of real weight $\alpha/2$ 
on the theta group with a unitary multiplier system. 
Hence $c_{\chi,\alpha}(m)=O(m^{\alpha/4})$ for $0< \alpha <4$. 
Therefore, we have 
\begin{equation}\label{328}
\aligned
Z_\alpha(s,\chi)
&= 2^{\alpha} \sum_{m=0}^{\infty} c_{\chi, \alpha}(m) \int_{0}^{\infty} e^{-\frac{\pi(m+\alpha) t^2}{q}} 
t^{s+\delta\alpha} \frac{dt}{t} \\
&= 2^{\alpha-1} q^{\frac{s+\delta\alpha}{2}} \, \pi^{-\frac{s+\delta\alpha}{2}}
\Gamma\bigl(\frac{s+\delta\alpha}{2}\bigr) 
\sum_{m=0}^{\infty} c_{\chi,\alpha}(m) \, (m+\alpha)^{-\frac{s+\delta\alpha}{2}}
\endaligned
\end{equation}
for ${\rm Re}(s)>(\alpha/2)+2-\delta \alpha$. 
From the functional equation $\eqref{lem2:02}$, 
the integrand $Z_\alpha(w,\chi) Z_{\beta}(s-w,\chi)$ on the right-hand side of $\eqref{332}$ 
is equal to $Z_\alpha(w,\chi)Z_\beta(w-s+\beta,\chi)$. 
From $\eqref{328}$, we obtain 
\begin{equation*}
\aligned
Z_\alpha & (w,\chi) Z_{\beta}(s-w,\chi) 
 = Z_\alpha(w,\chi)Z_\beta(w-s+\beta,\chi) \\
& =
\frac{1}{2} \, \left(\frac{\pi}{q}\right)^{\frac{s-\beta-\delta}{2}} 
\sum_{m=0}^{\infty} \sum_{n=0}^{\infty} 
c_{\chi,\alpha}(m) c_{\chi,\beta}(n) (m+\alpha)^{-\frac{\delta\alpha}{2}} (n+\beta)^{\frac{s-\beta(1+\delta)}{2}} \\ 
& \quad \times \Gamma\bigl(\frac{w+\delta\alpha}{2}\bigr) \Gamma\bigl(\frac{w-s+\beta+\delta\beta}{2}\bigr)
\, \bigr(q^{-1} \pi \sqrt{(m+\alpha)(n+\beta)}\,\bigr)^{-w}
\endaligned
\end{equation*}
for $w$ with a sufficiently large real part. 
Hence, if necessary, by moving the path of integration in $\eqref{332}$ to the right, 
we obtain 
\begin{equation}\label{333}
\aligned
\frac{1}{2\pi i} & \int_{c -i \infty}^{c + i \infty}
Z_{\alpha}(w,\chi)Z_\beta(s-w,\chi) dw 
= 2 \,\sum_{m=0}^{\infty} \sum_{n=0}^{\infty} 
c_{\chi,\alpha}(m) c_{\chi,\beta}(n) \bigl( \frac{n+\beta}{m+\alpha} \bigr)^{\frac{s-\beta+\delta(\alpha-\beta)}{4}} \\
& \qquad \qquad \qquad \qquad \qquad \qquad \qquad \quad 
\times K_{\frac{s-\beta+\delta(\alpha-\beta)}{2}}\bigl(\, \frac{2 \pi}{q} \sqrt{(m+\alpha)(n+\beta)} \,\,\bigr)
\endaligned
\end{equation}
by using $\eqref{lem5:02}$ of Lemma 5. 
We find that the series on the right-hand side of $\eqref{333}$ 
converges absolutely for any $s \in {\Bbb C}$ by a reason similar to that in the proof of Theorem 1. 
By combining $\eqref{332}$ with $\eqref{333}$, we obtain Theorem 2. \hfill $\Box$ 

\subsection{Proof of Theorem 3} 
From Lemma 3, we can apply Lemma 4 to 
$f(\frac{it}{\sqrt{q}})^\alpha$ and $f(\frac{it}{\sqrt{q}})^\beta$. 
By a way similar to that in the proof of Theorem 2, 
we obtain
\begin{equation}\label{338}
L^\ast(s,f)
=\frac{1}{2\pi i}\int_{c -i \infty}^{c + i \infty}
Z_{\alpha}(w, f)Z_\beta(s-w, f) dw
\end{equation}
for any positive $\alpha$ and $\beta$ with $\alpha+\beta$ and $s \in {\Bbb C}$, 
where $c$ is an arbitrary real number.  
Because $f(z)^\alpha$ is a modular form of real weight $k \alpha$ 
on $\Gamma_0(q)$, $c_{f,\alpha}(m)$ is at most of polynomial order (cf. ~\cite[\S3.2]{LaRa}).  
Hence, we have, by $\eqref{222}$, 
\begin{equation}\label{337}
\aligned
Z_\alpha(s,f)
&= q^{\frac{s}{2}} \, (2\pi)^{-s} \Gamma(s)
\sum_{m=0}^{\infty} c_{f,\alpha}(m) \, (m+\alpha)^{-s},
\endaligned
\end{equation}
when ${\rm Re}(s)$ is sufficiently large. 
From the functional equation $\eqref{lem3:02}$, 
$\eqref{337}$ and the assumption $Wf=f$, we have 
\begin{equation}\label{339}
\aligned
Z_\alpha(w,f) & Z_{\beta}(s-w,f)
= Z_\alpha(w,f)Z_\beta(w-s+k \beta,f) \\
= \, &
q^{-\frac{s-k \beta}{2}} 
(2\pi)^{s-k \beta}
\sum_{m=0}^{\infty} \sum_{n=0}^{\infty} 
c_{f,\alpha}(m) c_{f,\beta}(n) 
(n+\beta)^{s-k \beta} \\ 
& \times \Gamma(w) \Gamma(w-s+ k\beta)\,
\bigr( (2\pi)^2 \, q^{-1} (m+\alpha)(n+\beta) \,\bigr)^{-w}.
\endaligned
\end{equation}
By using $\eqref{lem5:03}$ of Lemma 5, we have
\begin{equation}\label{340}
\aligned
\frac{1}{2\pi i} & \int_{c -i \infty}^{c + i \infty}
Z_{\alpha}(w,f)Z_\beta(s-w,f) dw \\
=& 2 \,\sum_{m=0}^{\infty} \sum_{n=0}^{\infty} 
c_{f,\alpha}(m) c_{f,\beta}(n) \bigl( \frac{n+\beta}{m+\alpha} \bigr)^{\frac{s-k \beta}{2}} 
K_{s-k \beta}\bigl(\, \frac{4 \pi}{\sqrt{q}} \sqrt{(m+\alpha)(n+\beta)} \,\,\bigr).
\endaligned
\end{equation}
By combining $\eqref{338}$ with $\eqref{340}$, we obtain Theorem 3. \hfill $\Box$ 


\section{Proof of Theorem 4}
In this section, we prove Theorem 4 in a more general form. 
To state such a more general theorem, 
we introduce the function $E_{\alpha,\beta}(y,s)$ as follows.  

Let $\alpha$ and $\beta$ be positive real numbers. 
Define the function
$E_{\alpha,\beta}(y,s)$ by 
\begin{equation}\label{501}
E_{\alpha,\beta}(y,s)=
\int_{0}^{\infty} \bigl(\theta(y^2t^2)^\alpha \theta(y^{-2}t^2)^\beta -1 \bigr)\, t^s \frac{dt}{t}.
\end{equation}
The integral on the right-hand side converges absolutely for 
$y>0$ and $s$ with ${\rm Re}(s)>\alpha+\beta$, 
because the integrand decays exponentially fast as $t \to +\infty$ 
and is estimated as $O(t^{-\alpha-\beta})$ as $t \to +0$.  
From the definition, we find that $E_{\alpha,\beta}(1,s)=\zeta^\ast(s)$, 
whenever $\alpha+\beta=1$. In particular, $E_{1/2,1/2}(1,s)=\zeta^\ast(s)$.

\begin{thm}
Let $\alpha$ be positive real number. 
Denote $E_{\alpha,\alpha}(y,s)$ by $E_{\alpha}(y,s)$. 
Then,  
\begin{enumerate}
\item as a function of $y$, $E_{\alpha}(y,s)$ satisfies the modular equation $E_{\alpha}(y^{-1},s)=E_{\alpha}(y,s)$, 
\item as a function of $s$, $E_{\alpha}(y,s)$ is a meromorphic function on $\Bbb C$, 
\item as a function of $s$, $E_{\alpha}(y,s)$ is holomorphic except for the simple poles at $s=0$ and $s=2 \alpha$ 
with residues $-1$ and $1$, respectively, 
\item $E_{\alpha}(y,s)$ satisfies the functional equation $E_{\alpha}(y,s)=E_{\alpha}(y,2\alpha-s)$, 
\item $E_{\alpha}(y,s)$ has the expansion 
\begin{equation}\label{509}
\aligned
E_{\alpha}(y,s)
&= Z_\alpha(s) \, y^{s} + Z_\alpha(2\alpha-s) \, y^{2\alpha -s} \\
& \quad +y^{\alpha} \sum_{m=1}^{\infty} \sum_{n=1}^{\infty} 
c_\alpha(m) c_\alpha(n) \bigl( \frac{n}{m} \bigr)^{\frac{s-\alpha}{4}} 
K_{\frac{s-\alpha}{2}} \bigl( 2 \pi y^2 \sqrt{mn} \,\bigr)
\endaligned
\end{equation}
for any $y>0$ and $s \in {\Bbb C}$, 
where $Z_\alpha(s)$ and $c_\alpha$ are defined in $\eqref{106}$ and $\eqref{107}$, respectively, 
\item the singularities of $E_{\alpha}(y,s)$ depend only on 
the first two terms on the right-hand side of $\eqref{509}$. 
\end{enumerate}
\end{thm}
Clearly, Theorem 5 leads Theorem 4 by taking $\alpha=1/2$. 
The above properties of $E_{\alpha}(y,s)$ resemble the properties of $E^\ast(z,s)$. 
In this sense, we regard $E_{\alpha}(y,s)$ as an analogue of $E^\ast(z,s)$.  
From this point of view, 
we consider $\eqref{509}$ as an analogue of the Chowla-Selberg formula $\eqref{103}$. 
As a corollary of $\eqref{509}$, we obtain the following limit formula. 

\begin{cor}[an analogue of Kronecker's limit formula]
Let $\alpha$ be a positive real number. For any fixed $y>0$, we have
\begin{equation}
\aligned
\lim_{s \to 2 \alpha}  \Bigl(\, E_{\alpha}(y,s) - \frac{1}{s-2\alpha} \,\Bigr) 
&= Z_{\alpha}(2\alpha) \, y^{2\alpha} + \gamma_\alpha  \\ 
& + y^{\alpha} \sum_{m=1}^{\infty} \sum_{n=1}^{\infty} 
c_\alpha(m) c_\alpha(n) \bigl( \frac{n}{m} \bigr)^{\frac{\alpha}{4}} 
K_{\frac{\alpha}{2}} \bigl( 2 \pi y^2 \sqrt{mn} \,\bigr), 
\endaligned
\end{equation} 
where the constant $\gamma_\alpha$ is given by 
\begin{equation}
\gamma_\alpha=\lim_{s \to \alpha} \Bigl(\, Z_\alpha(s) -\frac{1}{s - \alpha} \, \Bigr).
\end{equation}
\end{cor}

\noindent
{\bf Proof of Theorem 5}. 
From definition $\eqref{501}$, we have 
\begin{equation}\label{502b}
E_{\alpha,\beta}(y^{-1},s) = E_{\beta,\alpha}(y,s). 
\end{equation}
Because $E_\alpha(y,s)=E_{\alpha,\alpha}(y,s)$, we obtain (1). 
By splitting integral $\eqref{501}$ 
into the part from $1$ to $\infty$ and the part from $0$ to $1$, 
we have
\begin{equation}\label{504}
\aligned
E_{\alpha,\beta}(y,s)
= & \, \frac{y^{\beta-\alpha}}{s-\alpha-\beta} - \frac{1}{s} 
 + \int_{1}^{\infty} \bigl( \theta(y^2t^2)^{\alpha} \theta(y^{-2}t^2)^{\beta} -1 \bigr) \, t^{s}\frac{dt}{t} \\
&+ y^{\beta-\alpha} 
\int_{1}^{\infty} \bigl( \theta(y^{-2} t^2)^{\alpha} \theta(y^{2}t^2)^{\beta} -1 \bigr) \, t^{\alpha+\beta-s}\frac{dt}{t}.
\endaligned
\end{equation}
Here we used the functional equation of $\theta(x)$. 
Representation $\eqref{504}$ gives a meromorphic continuation of $E_{\alpha,\beta}(y,s)$ 
and the functional equation 
\begin{equation}\label{505}
E_{\alpha,\beta}(y,\,s)=y^{\beta-\alpha} E_{\beta,\alpha}(y,\,\alpha+\beta-s). 
\end{equation}
By taking $\beta=\alpha$, we obtain (2), (3) and (4). 
To prove (5), we introduce
\begin{equation}\label{502}
F_{\alpha,\beta}(y,s)=
\int_{0}^{\infty} \bigl( \theta(y^2t^2)^\alpha -1 \bigr) \bigl( \theta(y^{-2}t^2)^\beta -1 \bigr) \,t^s \frac{dt}{t}.
\end{equation}
The integral on the right-hand side converges absolutely 
for $y>0$ and ${\rm Re}(s)>\alpha+\beta$ 
by a reason similar to that for $\eqref{501}$. 
We find that $E_{\alpha,\beta}(y,s)$ and $F_{\alpha,\beta}(y,s)$ are related as 
\begin{equation}\label{503}
F_{\alpha,\beta}(y,s) 
= E_{\alpha,\beta}(y,s)
 - Z_\alpha(s) \, y^{-s} - Z_\beta(s) \, y^{s}
\end{equation}
for ${\rm Re}(s)>\alpha+\beta$ from the identity similar to $\eqref{316}$. 
By a way similar to the proof of Theorem 1, we have 
\begin{equation}\label{508a}
\aligned
F_{\alpha,\beta}(y,s)
=& \, Z_\alpha(s-\beta) \, y^{2\beta -s} - Z_\alpha(s) \, y^{-s} 
+ \frac{1}{2\pi i} \int_{c-i \infty}^{c + i\infty} 
Z_\alpha(w) Z_{\beta}(s-w) y^{s-2w} dw \qquad 
\endaligned
\end{equation}
for ${\rm Re}(s)>\alpha+\beta$ and $c>{\rm Re}(s)$. 
Hence, we obtain
\begin{equation}\label{508b}
\aligned
E_{\alpha,\beta}(y,s)
= \, Z_\beta(s) \, y^{s} + Z_\alpha(s-\beta) \, y^{2\beta -s} 
+ \frac{1}{2\pi i} \int_{c-i \infty}^{c + i\infty} 
Z_\alpha(w) Z_{\beta}(s-w) y^{s-2w} dw \qquad 
\endaligned
\end{equation}
for ${\rm Re}(s) > \alpha+\beta$ and $c > {\rm Re}(s)$. 
Here, we find that
\begin{equation}\label{508c}
\aligned
\frac{1}{2\pi i} \int_{c-i \infty}^{c + i\infty} 
Z_\alpha(w) & Z_{\beta}(s-w) y^{s-2w} dw \\
&= y^{\beta} \sum_{m=1}^{\infty} \sum_{n=1}^{\infty} 
c_\alpha(m) c_\beta(n) \bigl( \frac{n}{m} \bigr)^{\frac{s-\beta}{4}} 
K_{\frac{s-\beta}{2}} \bigl( 2 \pi y^2 \sqrt{mn} \,\bigr)
\endaligned
\end{equation}
by an almost same way as in the proof of Theorem 1. 
By combining $\eqref{508b}$ with $\eqref{508c}$, we obtain
\begin{equation}\label{508}
\aligned
E_{\alpha,\beta}(y,s)
&= \, Z_\beta(s) \, y^{s} + Z_\alpha(s-\beta) \, y^{2\beta -s} \\
& \quad +y^{\beta} \sum_{m=1}^{\infty} \sum_{n=1}^{\infty} 
c_\alpha(m) c_\beta(n) \bigl( \frac{n}{m} \bigr)^{\frac{s-\beta}{4}} 
K_{\frac{s-\beta}{2}} \bigl( 2 \pi y^2 \sqrt{mn} \,\bigr)
\endaligned
\end{equation}
for any $y>0$ and $s$ with ${\rm Re}(s) > \alpha +\beta$. 
Because the series on the right-hand side converges absolutely for any $s \in {\Bbb C}$, 
$\eqref{508}$ holds for any $y>0$ and $s \in {\Bbb C}$ by analytic continuation.    
By taking $\beta=\alpha$, we obtain (5) and (6). 
\hfill $\Box$



\bigskip \noindent
Graduate School of Mathematics,\\
Nagoya University,\\
Chikusa-ku, Nagoya 464-8602,\\
Japan\\
e-mail address:m99009t@@math.nagoya-u.ac.jp

\vfill
\pagebreak

\begin{figure} \begin{center}
 \epsfile{file=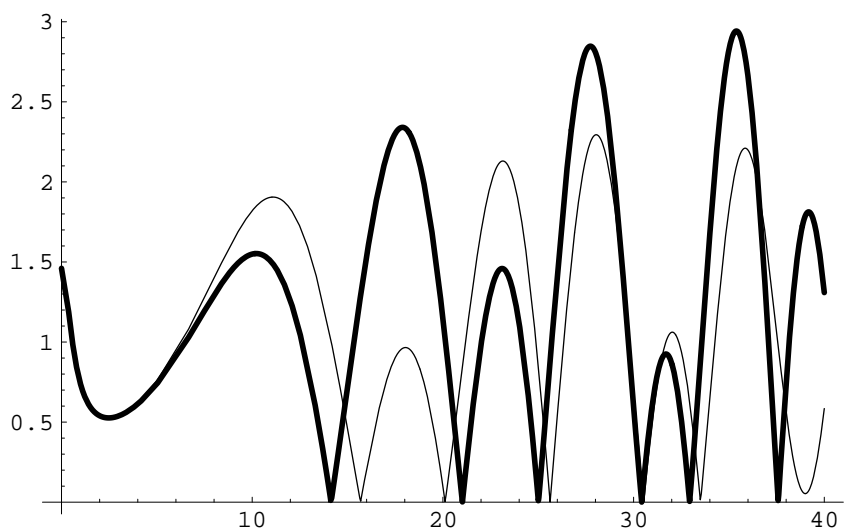,height=7cm}
 \caption{
 The thick line is $|\zeta(1/2+it)|$, and the thin line is 
 the absolute value of the sum of the first two terms 
 on the right-hand side of $\eqref{fig2}$ with 
 $s=1/2+it$ for $ 0 \leq t \leq 40$. } 
\end{center} \end{figure}

\begin{figure} \begin{center}
 \epsfile{file=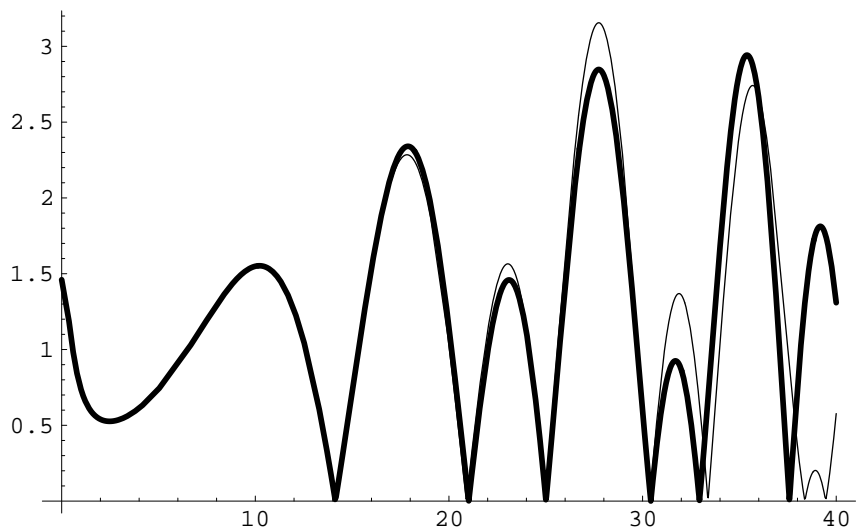,height=7cm}
 \caption{
 The thick line is $|\zeta(1/2+it)|$ and the thin line is 
 the absolute value 
 of the sum of the first two terms 
 and the terms corresponding to $1 \leq n \leq 3$ 
 of the infinite series 
 on the right-hand side of $\eqref{fig2}$ 
 with $s=1/2+it$ 
 for $0 \leq t \leq 40$. }
\end{center} \end{figure}

\end{document}